\documentclass[12pt]{article}
\usepackage{amssymb}
\usepackage{times,a4,amsmath,amsthm}
\usepackage[11pt,curve,matrix,arrow,frame,graph]{xy}

\newcommand\lnms{\emph{ Lecture Notes in Math. Springer-Verlag, Berlin, New York}}

\newcommand\cmp{\emph{ Comm. math. Phys.}}

\newcommand\jfa{\emph{ J. Funct. Anal.}}

\newcommand\jdg{\emph{ J. Differential Geom.}}

\def \dim{\mathop{\rm dim}\nolimits}
\def \Isom{\mathop{\rm Isom}\nolimits}

\def \ind{\mathop{\rm ind}\nolimits}
\newtheorem{rem}[equation]{Remark}
\newtheorem{theo}{Theorem}[section]
\newtheorem{pro}[theo]{Proposition}
\newtheorem{lem}[theo]{Lemma}
\newtheorem{cor}[theo]{Corollary}
\newtheorem{de}[theo]{Definition}

\def \reel{ {\rm I}\!{\rm R} }
\def \Sphere{ {\bf S} }

\def\ent{{{\rm Z}\mkern-5.5mu{\rm Z}}}
\def\vol{{\rm vol}\,}
\def\la{{\longrightarrow}}

\def \ima{{\rm Im}}

\def\no{${\rm n}^o$}

\def \cou {{\cal R}}
\def \Disk{{\bf D}}
\def \T {{\bf T}}
\def \H {{\mathcal H}}
\def \proof {\par{\sl Proof}\pointir}
\def  \sob{\forall u\in C^{\infty}_0(M),\ \,\mu _p(M)\left(\int_M|u|^{\frac{2p}{ p-2}}(x)dx\right)^{1-\frac{2}{
p}}\le\int_M|du|^2(x) dx,}
\def\ent{{{\rm Z}\mkern-5.5mu{\rm Z}}}
\def\proof{{\it Proof }.-- }
\def\endproof{\hfill {\bf Q.E.D.}}

\title{$L^2$ cohomology of Manifolds with flat ends.}

\author{Gilles Carron\\
D\'epartement de Math\'ematiques U.M.R. \\
 Universit\'e de Nantes \\
 2 rue de la Houssini\`ere \\
F-44322 Nantes cedex 03  \\
 Gilles.Carron@math.univ-nantes.fr \\
}

\begin{document}
\maketitle

\begin{abstract}
We give a topological interpretation of the space of $L^2$-harmonic forms on 
Manifold with flat ends. It is an answer to an old question of J. Dodziuk.
We also give a Chern-Gauss-Bonnet formula for the $L^2$-Euler characteristic
of some of these Manifolds. These results are applications of general theorems on
complete Riemannian Manifold whose Gauss-Bonnet operator is non-parabolic at infinity.
\end{abstract}

\centerline{\bf R\'esum\'e:}
\vskip0.25cm
Nous donnons une interpr\'etation topologique des espaces de formes harmoniques $L^2$ d'une
vari\'et\'e riemannienne compl\`ete plate à l'infini. Ceci r\'epond à une question pos\'ee par J. Dodziuk.
Nous obtenons aussi une formule de Chern-Gauss-Bonnet pour la caract\'eristique d'Euler $L^2$ de
certaines de ces vari\'et\'es. Ces r\'esultats sont des cons\'equences de th\'eor\`emes g\'en\'eraux sur les
vari\'et\'es riemanniennes compl\`etes dont l'op\'erateur de Gauss-Bonnet est non-parabolique à l'infini.

\vskip1cm
{\it Mathematics Subjet Classification (2000) : 58J10, 35J25, 53C21, 58C40.}

\section{Introduction}

Let $(M^n,g)$ be a complete Riemannian Manifold, we note $\H^k(M,g)$ or  $\H^k(M)$ its space of $L^2$-harmonic $k$-forms,
 that is to say the space of $L^2$ $k-$forms which are closed and coclosed. 
 These spaces have a (reduced) $L^2$-cohomology interpretation. When the Manifold is compact, without boundary,
  the celebrated theorem of Hodge-DeRham tells us that theses spaces are isomorphic to the real cohomology spaces of $M$ ;
   and we have the Chern-Gauss-Bonnet formula
$$\chi (M)=\sum_{k=0}^n(-1)^k\dim \H^k(M)=\int_M \Omega^g\, ,$$
where $\Omega^g$ is the Euler form of $(M^n,g)$ ; for instance if $\dim M=2$ then
$\Omega^g=\frac{KdA}{2\pi}$ , where $K$ is the Gaussian curvature and $dA$ the aera's form.

On non-compact Manifold, such results are generally not true ; for instance it is no longer true that the spaces
$\H^k(M,g)$ have finite dimension ; however it can happen it is the case. And the following questions are natural: 
\begin{enumerate}
\item What are the geometry insuring the finiteness of the dimension of the spaces $\H^k(M)$ ?
\item What are the links of the spaces $\H^k(M)$ 
with the topology of $M$ and with the geometry "at infinity" of $(M,g)$ ? 
\item And what kind of Chern-Gauss-Bonnet formula could we hope for the $L^2$- Euler characteristic
$$\chi_{L^2}(M)=\sum_{k=0}^n(-1)^k\dim \H^k(M)\ \ ?$$
\end{enumerate}
In (1982, \cite{Dod}), J. Dodziuk asked the following question: according to Vesentini (\cite{visentini})
 if $M$ is flat outside a compact set, the spaces $\H^k(M)$ are finite dimensional.
  Do they admit a topological interpretation ?

The main result of this paper solves this question. It's known that a complete Riemannian Manifold, which is flat outside
a compact set, has a finite number of ends. For sake of simplicity, we give, in the introduction, only the result
for Manifold with one flat end. 
\begin{theo} Let $(M^n,g)$ be a complete Riemannian Manifold with one flat end $E$ then 
\begin{enumerate}
\item If the volume growth of geodesic ball is at most quadratic :
$$\lim_{r\to\infty}\frac{\vol B_x(r)}{ r^2}<\infty$$ then we have 
$$\H^k(M,g)\simeq {\rm Im} \,\left(H^k_c(M)\longrightarrow H^k(M) \right).$$

\item If $\lim_{r\to\infty}\frac{\vol B_x(r) }{r^2}=\infty$, then the boundary of $E$
has a finite covering diffeomorphic to the product $\Sphere^{\nu-1}\times \T$ where $\T$ is a flat $(n-\nu)$-torus ; let
$\pi\,:\, \T\longrightarrow \partial E $ the induced immersion, then 
$$\H^k(M,g)\simeq H^k(M\setminus E, \ker \pi^*),$$
where  $H^k(M\setminus E, \ker \pi^*)$ is the cohomology associated to the complex of differential forms
on $M\setminus E$ which are zero when pull back to $\T$ : 
\begin{eqnarray*}
&H^k(M\setminus E, \ker \pi^*)=  \\ 
&\{\alpha \in C^\infty(\Lambda^k T^*(M\setminus E)),\ d\alpha=0,\ \pi^* \alpha=0\}
/
\{d\alpha,  \alpha\in C^\infty(\Lambda^{k-1} T^*(M\setminus E)),  \pi^* \alpha=0\}.
:\end{eqnarray*} 
\end{enumerate}

\label{the:hodge}\end{theo}

In potential theory, in the first case the manifold is parabolic, and in the second case it is non-parabolic (\cite{Ancona}).

This theorem was already known for asymptotically Euclidean Manifold, i.e. each end is simply connected 
(\cite{carrma, Melrose}).

This theorem is obtained through applications of the analysis we have developed in (\cite{carrpac,carrcre}) 
and with the help of the work of Eschenburg and Schroeder describing of the
 ends of such Manifold (\cite{ES}, see also \cite{ GPZ}).

Let's us explain our results on Dirac operator which are non-parabolic at infinity.
\begin{de} The Gauss-Bonnet operator $d+\delta$ of a complete Riemannian Manifold $(M,g)$ is
 called non-parabolic at infinity when
there is a compact set $K$ of $M$ such that for any bounded open subset  $U\subset M\setminus K$ there is a constant $C(U)>0$ with the inequality
\begin{equation}
\forall \alpha\in C^\infty_0(\Lambda T^*(M\setminus K)),\ \ C(U)\int_U|\alpha|^2\le \int_{M\setminus K}|d\alpha|^2+|\delta \alpha|^2.
\label{eq:weakf}
\end{equation}
\end{de}
The main property of these operators is the following
\begin{pro} If the Gauss-Bonnet operator of $(M,g)$ is non-parabolic at infinity then
$$\dim\{\alpha\in L^2(\Lambda T^*M),\ d\alpha=\delta\alpha=0\}<\infty.$$
Moreover let $D$ be a bounded open subset of $M$ containing $K$, and let $\,W(\Lambda T^*M)$
 be the Sobolev space which is the completion of
$C^\infty_0(\Lambda T^*M)$ with the norm 
$$\alpha\mapsto \int_D |\alpha|^2+\int_{M\setminus D} |d\alpha|^2+|\delta \alpha|^2,$$  
then this space imbedded continuously in $H^1_{loc}$ and
$$d+\delta\,:\,W(\Lambda T^*M)\la L^2(\Lambda T^*M)$$ is a Fredholm operator.
\end{pro}
 In fact, a differential form which is in the domain of $d+\delta$ is in $W$, that is to say
$$\{\alpha\in L^2, d\alpha\in L^2, \delta\alpha\in L^2\}\subset W.$$
So any $L^2$-harmonic forms is in $W$. The first step for proving our theorem \ref{the:hodge} is the following result

\begin{pro}If $(M^n,g)$ is a complete Riemannian Manifold whose curvatures vanish outside some compact set,
then the Gauss-Bonnet operator is non-parabolic at infinity.
\end{pro}

A special case which has been extensively studied is when the Gauss-Bonnet operator is Fredholm on its domain, or
equivalently when $0$ isn't in the essential spectrum of the Gauss-Bonnet operator. 
In this case according to N. Anghel (\cite{Anghel})
such operator is invertible at infinity, that is to say there exist a compact $K$ of $M$ and an constant $\Lambda>0$ such that
$$\forall \alpha\in C^\infty_0(\Lambda T^*(M\setminus K)),\ \ \Lambda\int_{M\setminus K}|\alpha|^2\le \int_{M\setminus K}|d\alpha|^2+|\delta \alpha|^2.$$
In this case, the Gauss-Bonnet operator is non-parabolic at infinity and the Sobolev space $W$ is the domain of $d+\delta$.
For instance, by the work of Borel and Casselman \cite{BC}, the Gauss-Bonnet operator 
is a Fredholm operator if $M$ is an even
dimensional locally symmetric space of finite volume and
negative curvature.

Another case is when the Manifold is with cylindrical end : that is to say
there is a compact $K$ of $M$ such that $M\setminus K$ is isometric to the Riemannian product $\partial K\times ]0,\infty[$.
According to the pioneering article of Atiyah-Patodi-Singer (\cite{APS}), the dimension of the
 space of $L^2$-harmonic forms is finite ; moreover these spaces are isomorphic to the image of the relative
cohomology in the absolute cohomology. These results were used by Atiyah-Patodi-Singer in order to obtain a 
formula for the signature
of compact Manifolds with boundary. In fact on Manifold with cylindrical end, the Gauss-Bonnet operator is 
non-parabolic at infinity; and what  Atiyah-Patodi-Singer have called a $L^2$ extended harmonic form is precisely a harmonic 
form which belongs to the Sobolev
space $W$. That's why we have called extended index the index of the operator
$$d+\delta\,:\,W(\Lambda T^*M)\la L^2(\Lambda T^*M).$$
In (\cite{carrcre}), we have developed analytical tools in order to compute this index. One of our results was that
this index only depends of the geometry of infinity. Recall that we have notice that
a harmonic $L^2$-form is in $W$ so that we have 
$$\begin{array}{lcl}
\ind_e(d+\delta)&=&\dim \{\alpha\in W(\Lambda T^*(M)),\, d\alpha+\delta\alpha=0\}\\
 & &-\dim \{\alpha\in L^2(\Lambda T^*(M)),\, d\alpha+\delta\alpha=0\}\\
 &=&\dim \frac{
\{\alpha\in W(\Lambda T^*(M)),\, d\alpha+\delta\alpha=0\}}
{\{\alpha\in L^2(\Lambda T^*(M)),\, d\alpha+\delta\alpha=0\} }\,.\hfill\\
\end{array}$$
In (\cite{carrcre}), we have shown the following 
\begin{theo} If $D$ is compact set such that outside we have the estimation
(\ref{eq:weakf}), let
 $$h_\infty(M\setminus D)=\dim \frac{\{\alpha\in W(\Lambda T^*(M\setminus D))\cap C^\infty,\, d\alpha+\delta\alpha=0\}
}{\{\alpha\in L^2(\Lambda T^*(M\setminus D))\cap C^\infty,\, d\alpha+\delta\alpha=0\}}$$
then we have 
$$h_\infty(M\setminus D)=2\ind_e (d+\delta).$$ 
\label{the:eind}
\end{theo}
A consequence of theorem on the topology of $M$ is the following exact sequence:

\begin{theo}If $(M,g)$ is a complete Riemannian Manifold whose Gauss-Bonnet operator is non-parabolic at infinity and
such that $\ind_e (d+\delta)=0$ that is to say such that every harmonic form in $W$ is in $L^2$ then
for every compact $D$ of $M$, we have the following exact sequence :
\begin{equation} 
...\la H^k(D,\partial D)\stackrel{\mbox{i}}{\longrightarrow} \H^k(M) \stackrel{\mbox{j}^*}{\longrightarrow} \H^k_{n}(M\setminus D)
\stackrel{\mbox{b}}{\longrightarrow} H^{k+1}(D,\partial D)\la...
\label{eq:suitexa}
\end{equation}
where $\H^k_{n}(M\setminus D)$ is the space of $L^2$ harmonic form on $M\setminus D$ whose normal components vanish along $\partial D$.
\label{the:thexa}
\end{theo}
 This exact sequence is the $L^2$-analogue of the exact sequence of the coboundary operator for the DeRham
cohomology. It's well known that this exact sequence is true for the non-reduced $L^2$-cohomology.
 This implies that if $0$ isn't in the essential spectrum of
the Gauss-Bonnet operator then the sequence (\ref{eq:suitexa}) holds, because in this case the (non-reduced)
 $L^2$-cohomology is isomorphic to the
space of $L^2$-harmonic forms ; note that in this case, the hypothesis of the theorem (\ref{the:thexa}) are satisfied,
 because $W$ is the domain of the Gauss-Bonnet operator.

When $(M,g)$ is a Manifold with one non-parabolic flat end, then the long exact sequence (\ref{eq:suitexa}) hold,
 and in this case the theorem (\ref{the:hodge}) is a consequence of the computation of the $L^2$-cohomology
  on the ends and of this exact sequence. 
This exact sequence can also be used to obtain a $L^2$-Chern-Gauss-Bonnet formula :
\begin{theo}  If $(M^n,g)$ is a complete oriented Riemannian Manifold of even dimension with one flat end $E$. Assume that
 $\lim_{r\to\infty}\frac{\vol B_x(r) }{r^2}=\infty$
 then 
$$\chi_{L^2}(M)=\int_M \Omega^g+q(E),$$
where $q(E)$ is computed in terms of $\pi_1(E)$ ; 
\begin{enumerate}
\item  When $\pi_1(E)$ has no torsion then $q(E)=0$.
\item  When ${\rm rank}\, \pi_1(E)=0$ we have $q(E)=\frac{1}{ |\pi_1(E)|}-1$.
\item  In general, $\pi_1(E)$ acts isometrically on the product $\Sphere^{\nu-1}\times \T$ where $\T$ is a flat $(n-\nu)$-torus 
and $\pi_1(E)\subset O(\nu)\times \left[\reel^{n-\nu}\rtimes O(n-\nu)\right]$, if $G_E$ is the image of $\pi_1(E)$ in
$O(n-\nu)$, then $$q(E)=-\frac{1}{|G_E|}\sum_{\gamma\in  G_E}\det({\rm Id}-\gamma).$$ 
\end{enumerate}
\end{theo}

This Gauss-Bonnet formula was already known for asymptotically Euclidean Manifold, i.e. each end is simply connected and
the curvatures almost vanish (cf. the work of Stern \cite{stern}, Borisov-M\"uller-Schrader \cite{BMS}, Br\"uning \cite{Br}) 
and also (\cite{carrduke}). 

This article is organized as follow, in a first part, we recall the main properties of the space of $L^2$-harmonics forms.
 In a second part, we present our analytical results of (\cite{carrpac,carrcre}) and we give examples. In the third part,
 we establish the long exact sequence (\ref{eq:suitexa}) and explain in which context this exact sequence is true.
The last part is devoted to the $L^2$-cohomology of Manifolds with flat ends.

{\flushleft\it Acknowledgements.} --- This paper was rewritten and finished when I was visiting
the Australian National University, at Canberra. It is a pleasure to thank this institution for its 
hospitality. I was supported in part by the "ACI" program of the French ministry of Research.

\section{The $L^2$-cohomology.}
We begin to recall what are the reduced $L^2$-cohomology spaces :
\subsection{Definition.} Let $(M^n,g)$ be a complete Riemannian Manifold of dimension $n$ : the operator of
  exterior differentiation is 
$$d\ :\ C^\infty_0(\Lambda ^kT^*M)\longrightarrow C^\infty_0(\Lambda ^{k+1}T^*M)$$ and it satisfies
$d\circ d=0$; its formal adjoint is
$\delta\ :\ C^\infty_0(\Lambda ^{k+1}T^*M)\longrightarrow C^\infty_0(\Lambda ^{k}T^*M)$; we have 
$$\forall \alpha \in C^\infty_0(\Lambda ^kT^*M),\ \forall\beta \in  C^\infty_0(\Lambda ^{k+1}T^*M),\ 
\int_M <d\alpha ,\beta >=\int_M <\alpha,\delta \beta> .$$
The spaces $Z^k_2(M)$ and $B^k_2(M)$ are defined as follow :
\begin{enumerate}
\item $Z^k_2(M)$ is the kernel of unbounded operator $d$ acting on 
$L^2(\Lambda ^kT^*M)$, or equivalently
$$Z^k_2(M)=\{\alpha \in L^2(\Lambda ^kT^*M),\ d\alpha =0\},$$ where the equation $d\alpha=0 $ has to
be understood in the distribution sense i.e. $\alpha\in Z^kL^2(M)$ if and only if
$$\forall\beta \in  C^\infty_0(\Lambda ^{k+1}T^*M),\ \int_M <\alpha,\delta \beta> =0\ .$$
That is to say $Z^k_2(M)= \left(\delta C^\infty_0(\Lambda ^{k+1}T^*M)\right)^\perp$.

\item $B^k_2(M)$ is the closure in $L^2(\Lambda ^kT^*M)$ of
$d\,\left[C^\infty_0(\Lambda^{k-1}T^*M)\right]$,  
where $C^\infty_0(\Lambda^{k-1}T^*M)$ is the space of smooth $k$-differential form with
compact support.

\end{enumerate}
 Because $d\circ d=0$, we have always $B^k_2(M)\subset Z^k_2(M)$, the $k-$
space of reduced $L^2$-cohomology is 
$$H^k_{2}(M)=Z^k_2(M)/B^k_2(M).$$ 
 Thus two weakly closed $L^2$ $k$-forms, $\alpha $ and $\beta $, are
$L^2$-cohomologuous  if and only if there is a sequence of smooth $(k-1)$-forms with compact
support, $\left(\gamma _l\right)_{l=0}^{\infty}$, such that $\alpha -\beta
=L^2-\lim_{l\to \infty}d\gamma _l.$

The  space of (non-reduced) $L^2-$cohomology is the quotient of
 $Z^k_2(M)$ by the following space
$$\{d\alpha,\ \alpha\in L^2(\Lambda^{k-1}T^*M),\ d\alpha\in L^2\}.$$
Now unless when it will be ambiguous, we will speak of $L^2$-cohomology for the reduced $L^2$-cohomology.

\subsection{ $L^2$-cohomology and harmonic forms.} 
 If $\H^k(M)$ is the space of $L^2$ harmonic $k$-forms :
$$\H^k(M)=\{\alpha \in L^2(\Lambda
^kT^*M),\ d\alpha =\delta \alpha =0\},$$ then the space $L^2(\Lambda ^kT^*M)$ has the following
of  Hodge-DeRham-Kodaira  orthogonal decomposition 
$$L^2(\Lambda ^kT^*M)=\H^k(M)\oplus\overline {d C^\infty_0(\Lambda ^{k-1}T^*M)}
 \oplus\overline {\delta C^\infty_0(\Lambda ^{k+1}T^*M)},$$
where the closure is taken with respect to  the  $L^2$ topology.
We also have 
$$Z^k_2(M)=\H^k(M)\oplus \overline {d C^\infty_0(\Lambda ^{k-1}T^*M)},$$
And so we have $$H^k_{2}(M)\simeq \H^k(M).$$

\subsection{Manifolds with compact boundary.}
If the Manifold is not complete and it has a compact boundary such that the Manifold,
together with its boundary, is metrically complete then we can also define absolute and relative $L^2$-cohomology 
and we have again an identification
of the $L^2$-cohomology space with a space of harmonic forms.  When the Manifold is compact with boundary 
these results are due to 
 G. Duff, D.C. Spencer \cite{DS}, and to P.E. Connor \cite{connor}, in fact the arguments of de G. Duff, D.C. Spencer
generalize easily to this more general case ; these results are well-known (see the works of
M. Lesch, J. Br\"uning \cite{BL}, J. Lott \cite{LottL2}, or \cite{carrduke}).

Let $(M^n,g)$ be a complete Riemannian Manifold
with compact boundary such that the Manifold, together with its boundary, is metrically complete.
We will note by $C^\infty_0(\Lambda^{k}T^*M)$  the space of smooth $k$-differential form with
compact support in the interior of $M$ and by $C^\infty_b(\Lambda^{k}T^*M)$ the space of smooth $k$-differential form with
bounded support in $M$, the elements of  $C^\infty_b(\Lambda^{k}T^*M)$ are not zero along the boundary.
We have two possible definitions for a $L^2$ form to be weakly closed or $L^2$-exact, we define the following four spaces:

\begin{enumerate}
\item $Z_{2,abs}^k(M)=\{\alpha \in L^2(\Lambda ^kT^*M),\  
\forall\beta \in  C^\infty_0(\Lambda ^{k+1}T^*M),\ \int_M <\alpha,\delta \beta> =0\ \}$
 or equivalently
$Z^k_{2,abs}(M)=\left(\delta C^\infty_0(\Lambda ^{k+1}T^*M)\right)^\perp$.

\item $B^k_{2,abs}(M)$ is the closure in $L^2(\Lambda ^kT^*M)$ of
$d\,\left[C^\infty_b(\Lambda^{k-1}T^*M)\right]$, 

\item $Z_{2,rel}^k(M)=\{\alpha \in L^2(\Lambda ^kT^*M),\  
\forall\beta \in  C^\infty_b(\Lambda ^{k+1}T^*M),\ \int_M <\alpha,\delta \beta> =0\ \}$
 or equivalently
$Z^k_{2,rel}(M)=\left(\delta C^\infty_b(\Lambda ^{k+1}T^*M)\right)^\perp$.

\item $B^k_{2,rel}(M)$ is the closure in $L^2(\Lambda ^kT^*M)$ of
$d\,\left[C^\infty_0(\Lambda^{k-1}T^*M)\right]$, 
\end{enumerate}
We have that a smooth form with bounded support is in $Z_{2,abs}^k(M)$ if and only if it is closed, and a
a smooth form with bounded support is in $Z_{2,rel}^k(M)$ if and only if it is closed and zero when pull-back
to the boundary. We have $B^k_{2,abs}(M)\subset Z^k_{2,abs}(M)$ and
$B^k_{2,rel}(M)\subset Z^k_{2,rel}(M)$, the absolute and relative (reduced) $L^2$-
cohomology spaces of $M$ are
$H^k_{2,abs}(M)=Z^k_{2,abs}(M)/B^k_{2,abs}(M)$ and $H^k_{2,rel}(M)=Z^k_{2,rel}(M)/B^k_{2,rel}(M)$.
Define $\H^k_{abs}(M)$ to be the space 
$$H^k_{2,abs}(M)=\{h\in L^2(\Lambda^{k}T^*M),\ dh=\delta h=0 {\rm\ and\ } {\rm int}_\nu h =0\},$$
where $\nu \,:\, \partial M\la TM$ is the inward unit normal field; and define
$\H^k_{rel}(M)$ to be the space
$$\H^k_{rel}(M)=\{h\in L^2(\Lambda^{k}T^*M),\ dh=\delta h=0 {\rm\ and\ } {\rm i}^* h =0\},$$
where $\ {\rm i}\,:\, \partial M\la M$ is the inclusion map.
We have then two orthogonal decompositions of $L^2(\Lambda^{k}T^*M)$ :

$$L^2(\Lambda^{k}T^*M)=\H^k_{abs}(M)\oplus\overline {d C^\infty_b(\Lambda ^{k-1}T^*M)}
 \oplus\overline {\delta C^\infty_0(\Lambda ^{k+1}T^*M)},$$
$$L^2(\Lambda^{k}T^*M)=\H^k_{rel}(M)\oplus\overline {d C^\infty_0(\Lambda ^{k-1}T^*M)}
 \oplus\overline {\delta C^\infty_b(\Lambda ^{k+1}T^*M)},$$
So we have the identifications
$$H^k_{2,abs}(M)\simeq\H^k_{abs}(M)\ {\rm and\ }H^k_{2,rel}(M)\simeq\H^k_{rel}(M).$$
where $\ {\rm i}\,:\, \partial M\la M$ is the inclusion map.

If we take the double of $M$ along $\partial M$ to obtain the Manifold $X=M\#_{\partial M} M$ with its
natural Lipschitz Riemannian metric, $X$ has a natural isometry the symmetry $\sigma$ which switches the two parts of
$M$ in $X$ ; and $\delta$ is well-defined so that we can speak of harmonic form on $(X,g\# g)$. 
Then 
\begin{pro}
The absolute $L^2$-cohomology of $(M,g)$ is naturally isomorphic to
the space of $\sigma$-invariant $L^2$-harmonic form on $(M\#_{\partial M} M, g\# g)$ and the relative
 $L^2$-cohomology of $(M,g)$ is naturally isomorphic to
the space of $\sigma$ anti-invariant $L^2$-harmonic form on $(M\#_{\partial M} M, g\# g)$.
\label{the:double}
\end{pro}
In fact, if $M^n$ is oriented then the Hodge star operator realizes an isomorphism between $H^k_{2,abs}(M)$
and $H^{n-k}_{2,rel}(M)$ and it is an isometry between $\H^k_{abs}(M)$
and $\H^{n-k}_{rel}(M)$.

When the Manifold has non-connected compact boundary, we can put the relative boundary condition
on some of the  connected compounents of the boundary and the absolute boundary condition on the others.

\subsection{The Bochner-Weitzenb\"ock formula}
When the Riemannian manifold $(M,g)$ is complete, and if $\Delta ^k=d\delta +\delta d$ is
the Hodge-DeRham Laplacian acting on  $k$-differential forms we have
$$\H^k(M)=\{\alpha \in L^2(\Lambda ^kT^*M),\ \Delta ^k \alpha =0\},$$
Moreover in any cases, complete or not we have the Bochner-Weitzenb\"ock decomposition
\begin{equation}
\Delta ^k=\bar\Delta +\cou^k,
\label{eq:BW}
\end{equation}
where $\cou^k$ is a symmetric linear operator of $\Lambda ^kT^*M$ which is defined with
the curvature operator of $(M^n,g)$ (cf \cite{GM}) ; for instance
$\cou^1$ is the  Ricci operator), and we have always the bound
$|\cou^k|(x)\le c(n)|R|(x)$ where $R$ is the curvature tensor of $(M,g)$ ; better if $\rho$ is the
curvature operator we have always $|\cou^k|(x)\le k(n-k)|\rho|(x)$.

\section{When does the $L^2$-cohomology space has finite dimension.}
The purpose of this part is to give conditions which insure that the space of harmonic $L^2$-forms
has finite dimension. These conditions are defined in \cite{carrpac}. In this paper, we study 
Dirac type operator which are non-parabolic at infinity. We give here the definition for the particular case of the Gauss-Bonnet operator.
\begin{de} The Gauss-Bonnet operator $d+\delta$ of a complete Riemannian Manifold $(M,g)$ is called non-parabolic at infinity  when
there is a compact set $K$ of $M$ such that for any bounded open subset $U\subset M\setminus K$ there is a constant $C(U)>0$ with the inequality
\begin{equation}
\forall \alpha\in C^\infty_0(\Lambda T^*(M\setminus K)),\ \ C(U)\int_U|\alpha|^2\le \int_{M\setminus K}|d\alpha|^2+|\delta \alpha|^2.
\label{eq:weakf2}
\end{equation}
\end{de}
\subsection{Properties of these operators.} 
In fact, we have the following
\begin{pro}
 If the Gauss-Bonnet operator of $(M,g)$ is non-parabolic at infinity then
$$\dim\{\alpha\in L^2(\Lambda T^*M),\ d\alpha=\delta\alpha=0\}<\infty.$$
\end{pro}
This has been proved in \cite{carrpac}, in this paper we have also give some characterizing properties
 of non-parabolicity at infinity, and in fact
 these operators satisfy the following
\begin{theo} If the Gauss-Bonnet operator of $(M,g)$ is non-parabolic at infinity, let $D$ be a bounded open subset of $M$ containing
the compact set $K$ outside which the estimates (\ref{eq:weakf2}) holds
and define $W(\Lambda T^*M)$ to be the Sobolev space which is the completion of
$C^\infty_0(\Lambda T^*M)$ with the quadratic form 
$$\alpha\mapsto \int_D |\alpha|^2+\int_{M\setminus D} |d\alpha|^2+|\delta \alpha|^2,$$  
then this space imbedded continuously in $H^1_{loc}$ and
$$d+\delta\,:\,W(\Lambda T^*M)\la L^2(\Lambda T^*M)$$ is Fredholm.
\end{theo}
The direct consequence of the Fredholmness of this operator is the following Hodge decomposition
\begin{equation}L^2(\Lambda^k T^*M)=\H^k(M)\oplus dW(\Lambda^{k-1}T^*M)\oplus \delta W(\Lambda ^{k+1}T^*M).\label{eq:hodgeW}
\end{equation}
And $Z^k_2(M)=\H^k(M)\oplus dW(\Lambda^{k-1}T^*M)$. Hence, we have
\begin{pro}
A $L^2$ closed form $\alpha$ is $L^2$ cohomologous
to zero if and only if there is a $\beta\in W(\Lambda^{k-1}T^*M)$ such that $\alpha=d\beta$ and we can always choose $\beta$ so that
 $\delta \beta=0$.
\label{pro:green}
\end{pro}

We have also the same feature for complete Riemannian Manifold with compact boundary for which the Gauss-Bonnet operator
is non-parabolic at infinity. Let $(M,g)$ be such a Riemannian Manifold, we define the space $W(\Lambda T^*M)$ 
has the completion of the space $C^\infty_b(\Lambda T^*M)$ with the respect to quadratic forms
$$\alpha\mapsto \|\alpha\|^2_{H^{\frac{1}{2}}(\partial M)}  +\|\alpha\|^2_{L^2(D)}+\|(d+\delta)\alpha\|^2_{L^2(M)};$$
where $D$ is a bounded open subset of $M$ containing the boundary of $M$ and containing the compact $K$
outside which we have the estimate (\ref{eq:weakf2}).
And we define $W_0(\Lambda T^*M)$ to be the closure of $C^\infty_0(\Lambda T^*M)$ in $W(\Lambda T^*M)$.
Then we have the orthogonal decompositions
$$L^2(\Lambda^k T^*M)=\H^k_{abs}(M)\oplus dW(\Lambda^{k-1}T^*M)\oplus \delta W_0(\Lambda ^{k+1}T^*M).$$

$$L^2(\Lambda^k T^*M)=\H^k_{rel}(M)\oplus dW_0(\Lambda^{k-1}T^*M)\oplus \delta W(\Lambda ^{k+1}T^*M).$$
\begin{pro}
 $\alpha\in Z^k_{2,abs}(M)$ is $L^2$ cohomologous
to zero if and only if there is a $\eta\in W(\Lambda^{k-1}T^*M)$ such that $\alpha=d\eta$ and we can always choose $\eta$ so that $\delta \eta=0$.

$\alpha\in Z^k_{2,rel}(M)$ is $L^2$ cohomologous
to zero if and only if there is a $\eta\in W_0(\Lambda^{k-1}T^*M)$ such that $\alpha=d\eta$ and we can always choose $\eta$ so that
  $\delta \eta=0$.
\label{pro:greenb}
\end{pro}
We note that if $(M,g)$ is a complete Riemannian Manifold for which the Gauss-Bonnet operator
is non-parabolic at infinity, then for every open set $U\subset M$ with smooth compact boundary, then
the Gauss-Bonnet operator on $(U,g)$ is non-parabolic at infinity; and moreover we have an exact sequence
$$0 \longrightarrow W_0(\Lambda T^*U)
  \longrightarrow W(\Lambda T^*M)\longrightarrow W(\Lambda T^*(M\setminus U))\longrightarrow 0.$$
  Where the first map is the extension by zero map and the second is the restriction map.
 
 The condition of being non-parabolic at infinity  depends only of the geometry in a neighborhood of infinity, so our result
of finiteness for the dimension of the space of $L^2$-harmonic form is in concordance with J. Lott's result [L2], which
asserts that the finiteness for the dimension of the space of $L^2$-harmonic form
 depends only of the geometry in a neighborhood of infinity : that is to say the spaces of $L^2$-harmonic form of
 two complete Riemannian Manifolds, which are isometric outside 
some compact set, have simultaneously finite or infinite dimension.

 Now we can give examples.
\subsection{Manifolds with flat ends.} In \cite{carrpac}, we have shown the following
\begin{pro}
 If $(M^n,g)$ is a complete Riemannian Manifold whose curvatures vanish outside some compact set,
then the Gauss-Bonnet operator is non-parabolic at infinity.
\end{pro}
 For sake of completeness, we recall the proof
 
\proof Let $K$ be a bounded open subset of $M$ outside which the curvatures vanish, then by the Bochner-Weitzenb\"ock
 formula we have
$$\forall \alpha\in C^\infty_0(\Lambda T^*(M\setminus K)), \ \ \int_{M\setminus K} |d\alpha|^2+|\delta\alpha|^2=\int_{M\setminus K} |\nabla \alpha|^2;$$
but according to the Kato lemma, we have $|\nabla\alpha|\ge |d|\alpha||$ so that we obtain
$$ \forall \alpha\in C^\infty_0(\Lambda T^*(M\setminus K)),\ \ \int_{M\setminus K} |d\alpha|^2+|\delta\alpha|^2\ge\int_{M\setminus K} \left|d| \alpha|\right|^2 .$$
Now if we consider the operator $H=\Delta+{\bf 1}_K$ acting on functions, $H$ is a positive operator and it satisfies the
inequality $\forall u\in C^\infty_0(M),\ \ \int_M <Hu,u>\ \ge \int_K |u|^2,$
then by a result on Ancona \cite{Ancona}, we know that for any bounded open subset $U$ of $M$ we have a positive constant $C(U)$ such that
$$\forall u\in C^\infty_0(M),\ \ \int_M <Hu,u>\ \ge C(U) \int_U |u|^2.$$
We take a bounded open subset $U$ of $M\setminus K$, and $\alpha\in C^\infty_0(\Lambda T^*(M\setminus K))$, we apply this inequality to $u=|\alpha|$,
 we have
$$\int_{M\setminus K} |d\alpha|^2+|\delta\alpha|^2\ge\int_{M\setminus K} \left|d| \alpha|\right|^2=
\int_{M} \left|d| \alpha|\right|^2+{\bf 1}_K|\alpha|^2\ge C(U)\int_U |\alpha|^2.$$
That is to say the Gauss-Bonnet operator is non-parabolic at infinity.\hfill\endproof

In fact we have only used the fact that the curvature operator is non negative outside some compact.
 
\begin{rem}
 If $D$ is a bounded open set containing $K$ then, a by-product of the proof is that an equivalent norm on $W$ is given by
 $$\alpha\mapsto \sqrt{ \int_{M} |\nabla \alpha|^2+\int_{D} | \alpha|^2}.$$
 This comes from the Bochner-Weitzenbock formula and standard elliptic estimate.
 Moreover if the Manifold $(M,g)$ is non-parabolic, then this norm   is equivalent to 
 $\alpha\mapsto \sqrt{ \int_{M} |\nabla \alpha|^2}$. As a matter of fact, a Riemannian Manifold is called
 non-parabolic if its Brownian motion is transient or equivalently, if it carries positive Green 
 functions; a analytical characterization has been given by A. Ancona (\cite{Ancona}):
 $(M,g)$ is non-parabolic if and only if for any bounded open subset  $U\subset M$ there is a positive constant
 $C(U)$ such that 
 $$ \forall f\in C^\infty_0(M)\ ,\ C(U) \int_U f^2\le \int_M |df|^2.$$
 This result is one of our source of inspiration to define non-parabolicity at infinity.
 So if $(M,g)$ is non-parabolic, the result of Ancona and the Kato inequality show that
 $W(\Lambda T^*M)$ is $H^1_0(\Lambda T^*M)$ the completion of $C^\infty_0(\Lambda T^*M)$
with the norm $\alpha\mapsto \|\nabla\alpha\|_{L^2}$ and we have that
the operator $d+\delta \,:\,H^1_0(\Lambda T^*M)\la L^2(\Lambda T^*M)$ is a Fredholm operator.
\label{rem:h10}
\end{rem}

 We can improve this result when the negative part of the curvature operator is controlled by the quadratic form
$\alpha\mapsto \|\nabla\alpha\|_{L^2}$. Let $\lambda(x)$ be the lowest eigenvalue of the operator $\cou(x)$ appearing in the
Bochner-Weitzenb\"ock formula (\ref{eq:BW}), we note by $\cou_-(x)$ the negative part of $\lambda(x)$ that is to say
$$
\cou_-(x)=\left\{\begin{array}{lcl } -\lambda(x) &\ {\rm if }\ & \lambda(x)\le 0 \\
                  0        &\ {\rm if }\ & \lambda(x)>0.\\ \end{array}\right.
                $$
\begin{pro} If $(M,g)$ is a complete Riemannian Manifold and if there is a compact set
$K$ of $M$ and a $\varepsilon>0$ such that
$$\forall \alpha\in C^\infty_0(\Lambda T^*(M\setminus K)),\ \ (1+\varepsilon)
\int_{M\setminus K} \cou_-(x)|\alpha|^2(x)dx\le \int_{M\setminus K}|\nabla\alpha|^2 $$
then the Gauss-Bonnet operator is non-parabolic at infinity. Moreover if $(M,g)$ is non-parabolic then
$d+\delta\,:\, H^1_0(\Lambda T^*M)\la L^2(\Lambda T^*M)$ is Fredholm.
\label{the:general}
\end{pro}
\proof If $ \alpha\in C^\infty_0(\Lambda T^*(M\setminus K))$ we apply the Bochner-Weitzenb\"ock formula
\begin{eqnarray*}
\int_{M\setminus K} |d\alpha|^2+|\delta\alpha|^2&=\int_{M\setminus K}|\nabla\alpha|^2+<\cou\alpha,\alpha>\\
                                          &\ge \int_{M\setminus K}|\nabla\alpha|^2-\cou_-(x)|\alpha |^2\\
                                          &\ge \frac{\varepsilon}{ 1+\varepsilon}\int_{M\setminus K}|\nabla\alpha|^2.
                                          \end{eqnarray*}
Then we have only to apply the argument given in the proof of the last proposition.
Moreover when the Riemannian Manifold is non-parabolic, the arguments of the remark (\ref{rem:h10}) show that
$H^1_0=W$.\hfill\endproof

\subsection{Sobolev-Orlicz inequalities.} One application of the proposition (\ref{the:general}) is based
 on the Sobolev-Orlicz inequality obtain in (\cite{carrcol}). Let's us recall 
what is this inequality. Let $(M,g)$ be an complete Riemannian Manifold and
let $\left(P(t,x,y)\right)_{(t\in \reel_+,\,x,y\in M)}$ be its heat kernel, that is to say this is the minimal 
solution of the Cauchy problem
$$\left\{
\begin{array}{cc}
\frac{\partial P}{ \partial t}(t,x,y)+\Delta^g_y P(t,x,y)=0, & \ (x,y)\in M,\, t>0 \\
P(0,x,y)=\delta_x(y).& \\
\end{array}\right.$$
Then let $\varphi$ be the function on $\reel^+\times M$ define by
$$\varphi(\lambda,x)=\lambda\left(\int^\infty_{\frac{1}{ 4\lambda}}\sqrt{\frac{P(s,x,x)}{ s}}ds\right)^2,$$
We assume that this integral is finite for some $x$ in $M$ (hence it's finite for all $x$ in $M$ 
by the Harnack inequality). 
Now if  $u\in C^\infty_0(M)$, we associate
$$N(u)=\sup\left\{\int_M uv,\ v\in C^\infty_0(M){\rm\ with\ } \int_M \varphi(|v|(x),x)dx \le 1\right\}\ ;$$
then $N$ is a norm and the completion of $C^\infty_0(M)$ with this norm is a Banach space (called an Orlicz space)
which is made of locally integrable functions ; moreover in \cite{carrcol}, we have shown the following universal Sobolev inequality
 $$\forall u\in C^\infty_0(M),\ \ N(u^2)\le C\|du\|_{L^2}^2,$$
this for an universal constant $C$. This shows that in this cases the Manifold is non-parabolic. As an application of this result we can state the following result
\begin{pro}
There is an universal constant $C$ so that if $(M,g)$ is a complete Riemannian Manifold, whose heat kernel
 $\left(P(t,x,y)\right)_{(t\in \reel_+,\,x,y\in M)}$ 
 and whose curvature operator $\cou$ satisfy 
$$\int_M \,\cou_-(x)\left(\int^\infty_ {\frac{C}{\cou_-(x)}}\sqrt{\frac{P(s,x,x)}{ s}}ds\right)^2 dx<\infty,$$
then $d+\delta\,:\,H^1_0(E)\la L^2(E)$ is Fredholm.
\label{the:orlicznpy}
\end{pro}

\proof By definition, we have the H\"older inequality
$$\int_Mu^2v\le N(u^2)\ \inf\left\{\lambda>0,\ \int_M\varphi\left(\frac{|v(x)|}{ \lambda},x\right)dx\le 1\right\}.$$
If $\int_M \varphi(2C\cou_-(x),x)dx<\infty$ then there is a compact $K$ of $M$ such that
$$\int_{M\setminus K} \varphi(2C\cou_-(x),x)dx\le 1.$$ Then if $\alpha \in C^\infty_0(\Lambda T^*(M\setminus K))$, we have
\begin{eqnarray*}<\cou_-(x)\alpha,\alpha>&\le \frac{1}{2C}N(|\alpha|^2)\\
&\le \frac{1}{ 2} \int_{M\setminus K}|\nabla\alpha|^2.
\end{eqnarray*}
this shows that $d+\delta\,:\,H^1_0(E)\la L^2(E)$ is Fredholm.\hfill\endproof

 For instance, if the Manifold satisfies the Sobolev inequality
$$\sob $$ then by the result of J. Nash [N], we know that the heat kernel satisfies a uniform bound
$$P(t,x,x)\le C/t^{p/ 2},\ \forall x\in M,\ \forall t>0$$
so that in this case, we have the uniform bound 
$$\varphi(\lambda,x)\le C \lambda^{p/2},\forall x\in M, \forall \lambda>0\ .$$
Thus the hypothesis of the theorem are satisfies if the curvature of $(M,g)$ is in $L^{p/ 2}$; 
so the proposition (\ref{the:orlicznpy}) generalizes the result of \cite{carrma}.

\subsection{The Warped product cases} We recall the results we have obtain in \cite{carrcre} for the Gauss-Bonnet operator on a Manifold
which is a warped product at infinity :
\begin{pro}
If $(M,g)$ is a complete Riemannian Manifold and if there is a compact set $K$ of $M$ such that
$(M\setminus K,g)$ is isometric to the warped product $\left(]0,\infty[\times\partial K,dr^2+f^2(r)g\right)$, then the Gauss-Bonnet operator is non-parabolic at infinity 
in the following two cases
\begin{enumerate}
\item$f(r)=ar$ for an $a>0$,
\item $\lim_{r\to \infty}f'(r)=0$.
\end{enumerate}
\label{the:wp}\end{pro}

\section{Topology and $L^2$-harmonic forms.}

\subsection{The exact sequence.}
The purpose of this paragraph is to give a general condition which insures that the
following sequence is exact : if $D$ is a compact domain of the complete Riemannian Manifold
$(M^n,g)$, we consider
\begin{equation}
..\la H^k(D,\partial D)\stackrel{\mbox{i}}{\longrightarrow} \H^k(M)\stackrel{\mbox{j}^*}{\longrightarrow} H^k_{2,abs}(M\setminus D)
\stackrel{\mbox{b}}{\longrightarrow} H^{k+1}(D,\partial D)\la..
\label{eq:suitexa2}
\end{equation}
This is the long sequence associated to the coboundary operator $\rm b$.  
First, let's us recall how the morphism ${\rm\, i\,}, j^*$ and $\rm b$ are defined :
\begin{enumerate}
\item ${\rm\, i\,} $ is the following natural application : 
to a closed smooth form with compact support in $D$, ${\rm\, i\,}$ associated its $L^2$-
cohomology class, that is to say
$${\rm\, i\,} \,\alpha\,{\rm mod } C^\infty_0(\Lambda^kT^*D)=\alpha\,{\rm mod}\overline{C^\infty_0(\Lambda^kT^*M)}^{L^2}.$$
\item $j^*$ is the morphism induced by the restriction map $j\,:\, M\la M\setminus D$.
\item The coboundary operator $\rm b$ is defined as usually :
for each cohomology class of $M\setminus D$, there is a smooth form $\alpha$ in it which is closed on an neighborhood of $\partial D$, 
${\rm b}[\alpha]$ is the cohomology class of $d\bar \alpha$ where
$\bar \alpha$ is a smooth extension of $\alpha$ which is closed in a neighborhood of $\partial D$. $\rm b$ is well defined,
 it does not 
depend on the choice of $\alpha$ and nor of its extension.
\end{enumerate}
 By construction, we have $${\rm j}^*\circ{\rm  i}=0,\  {\rm b}\circ{\rm  j}^*=0 {\rm\ and\ }
 {\rm i}\circ{\rm  b}=0\ ;$$ so we have the inclusions
$${\rm Im\ i}\subset{\rm Ker\ j}^*,\ {\rm Im\ j}^*\subset{\rm Ker\ b}{\rm\ et\  }{\rm Im\ b}\subset{\rm Ker\ i}.$$
  In \cite{carrduke}, we have noticed that
\begin{pro}
The equality $\ker {\rm b}=\ima j^*$ always hold.
\end{pro}
\proof This comes from the long exact sequence in the DeRham cohomology. If $[\beta]\in H^k_2(M\setminus D)$ is in the kernel of $\rm b$
there is a smooth extension of $\beta$ say $\bar \beta$, and a smooth (k-1)-form $\gamma$ with compact support in $D$ so that
$d\bar\beta=d\gamma$ ; now the form $\bar\beta-\gamma$ is closed, square integrable and 
his restriction to $M\setminus D$ is $\beta$ so that
$[\beta]=j^*[\bar\beta-\gamma]$.\hfill\endproof

In fact when the Gauss-Bonnet operator is non-parabolic at infinity  we have more
\begin{pro}If the Gauss-Bonnet operator is non-parabolic at infinity , then 
$$\ker j^*=\ima {\rm\, i\,}.$$
\end{pro}
\proof Let $h$ be a harmonic $L^2$ k-form on $M$ such that its restriction to $M\setminus D$ is $L^2$-cohomologous to zero.
 By (\ref{pro:greenb}), we have a smooth form $\eta\in W(\Lambda^{k-1}T^*(M\setminus D))$ 
 with $h=d\eta$ and $\delta \eta=0$ on $M\setminus D$.
Now if $\overline{\eta}\in W(\Lambda^{k-1}T^* M)$ is an smooth extension of  $\eta$
then $h-d\overline{\eta}$ is a closed $L^2$ form which is $L^2$ cohomologous to $h$.
Then the support of  $h-d\overline{\eta}$ is in $\overline{D}$. And it is zero when pull-back to $\partial D$.
Hence the $L^2$ cohomology class of $h$ is in the image of the natural map ${\rm\, i\,}$ from
$H_c^k(D)\simeq H^k(D,\partial D)$ in $\H^k(M)$ .\hfill\endproof

For the remaining equality, we conclude with the analysis done in \cite{carrcre} :
\begin{theo}
If the Gauss-Bonnet operator is non-parabolic at infinity and if any harmonic form in $W$ is in $L^2$ then 
the long sequence (\ref{eq:suitexa2}) is exact.
\end{theo}
\proof Our hypothesis is
$$h_\infty(M)=\dim
\frac
{ \{\alpha\in W,\, d\alpha+\delta\alpha=0\} }
{ \{\alpha\in L^2,\, d\alpha+\delta\alpha=0\} }
=0.$$
In \cite{carrcre}, we have computed these dimension, our result is that if we define
$$h_\infty(M\setminus D)=\dim
\frac
{\{\alpha\in W(\Lambda T^*(M\setminus D))\cap C^\infty,\, d\alpha+\delta\alpha=0\}}
{\{\alpha\in L^2(\Lambda T^*(M\setminus D))\cap C^\infty,\, d\alpha+\delta\alpha=0\}}
$$ then we have
$$h_\infty(M)=\frac{1}{ 2}h_\infty(M\setminus D).$$
Our hypothesis implies that a smooth harmonic form on $M\setminus D$ which is in $W$ is in fact squared integrable. We can now show that 
$\ker {\rm\, i\,}\subset\ima {\rm b}$. 
Let $\alpha$ be a closed $k$-form with compact support in $D$ which is zero in $L^2$-cohomology, thus by 
(\ref{pro:greenb}), we have 
a $\beta\in W(\Lambda^{k-1}T^*M)$ with
$$\delta \beta=0\ {\rm and}\ \alpha=d\beta$$
Now $\beta$ is smooth by elliptic regularity, and the restriction of $\beta $ to $M\setminus D$ is harmonic and is in $W$;
so that $\beta$ is in fact squared integrable by our hypothesis.
and hence $[\alpha]={\rm b}[\beta|_{M\setminus D}]$.\hfill\endproof 

In fact, we can also consider this long exact sequence if $(M,g)$ is a complete Riemannian Manifold with compact boundary
$\partial M$ and if we put on the boundary appropriate boundary condition (the relative condition
on some connected component of $\partial M$ and the absolute condition on the other), and if
$U\subset M$ is a open set with smooth compact boundary with $\partial U\subset {\rm int }(M)$, then we can define
the sequence 
$$..\la H^k_{2,rel}(U) \stackrel{\mbox{i}}{\longrightarrow} H^k_{2}(M)
 \stackrel{\mbox{j}^*}{\longrightarrow} H^k_{2,abs}(M\setminus U)
\stackrel{\mbox{b}}{\longrightarrow} H^{k+1}_{2,rel}(U)\la..$$
On $M\setminus U$,  $H^k_{2,abs}(M\setminus U)$ is defined with the absolute boundary condition on
$\partial U$ and on $\partial M$ we put the same boundary condition as for defining  $H^k_{2}(M)$.
The map $i$ is well defined because by definition an element of $Z^k_{2,rel}(U)$ when extended by zero on $M$
is an element of $Z^k_{2}(M)$. For the coboundary map ${\rm b}$, due to the non compactness of $U$, we have to consider
only extension with compact support, in fact this coboundary map is the
composition of the natural map from $H^k_{2, abs}(M\setminus U)$ to the DeRham cohomology group $H^k_{ abs}(M\setminus U)$ then of the
usual coboundary map from $H^k_{ abs}(M)$ 
to $H^{k+1}_{ c}(U)$ (the group of cohomology with compact support of $U$) and finally of the natural map from
$H^{k+1}_{ c}(U)$ to $H^{k+1}_{2,rel}(U)$.
The same proof leads to the following results

\begin{theo}
If the Gauss-Bonnet operator is non-parabolic at infinity then $\ker {\rm j}^*=\ima\, {\rm i}$. Moreover if
 if any harmonic form in $W(\Lambda T^*M)$ is in $L^2(\Lambda T^*M)$ then 
the long exact sequence is true.
More exactly if $h_\infty(U)=0$ then $\ker {\rm b}=\ima j^*$ and if $h_\infty(M\setminus U)=0$ then 
$\ker {\rm i}=\ima {\rm b}$.
\label{the:seb}
\end{theo}

The fact that  $\ker {\rm j}^*=\ima \,{\rm i}$ implies the following
\begin{cor}If the Gauss-Bonnet operator is non-parabolic at infinity and if $U$ is an open set of $M$ with compact boundary then
$\ima\left( H^k_{2,rel}(U) \longrightarrow H^k_{2,abs}(U)\right)$
injects in $H^k_{2}(M)$.
\label{the:relabs}
\end{cor} 
\proof If $[\alpha]\in H^k_{2,rel}(U)$ is map to zero in $H^k_{2}(M)$, then by (\ref{pro:greenb}) there is a $\beta\in W(\Lambda^{k-1} T^*M)$, such that
$\alpha=d\beta$ on $M$, this identity is also true on $U$, and
$\beta|_U\in W(\Lambda^{k-1} T^*U)$. So $[\alpha]$ is map to zero in $H^k_{2,abs}(U)$
\endproof

If $U=M$, then the result is due to Anderson (\cite{Anderson}): 
 $\ima\left(  H^k_{rel}(M) \longrightarrow H^k_{abs}(M)\right)$
always injects in $H^k_{2}(M)$.

\subsection{Some examples.}
In \cite{carrduke}, we have shown the following 
\begin{theo}
 If $(M^n,g)$ is a complete Riemannian Manifold, which for a $p>4$ satisfies the Sobolev inequality
$$\sob,$$
and whose curvature tensor $R$ satisfies 
$$\int_M|R|^{p/ 2}(x)dx<\infty$$
then the sequence (\ref{eq:suitexa2}) is exact.
\end{theo}
This result can be recover by the analysis done here, because as we notice earlier,
 on such Manifold the Gauss-Bonnet operator
is non-parabolic at infinity. Moreover, the main analytical tool used in \cite{carrduke} was that if $\alpha\in H^1_0$ satisfies
$\Delta\alpha=(d\delta+\delta d)\alpha\in C^\infty_0(\Lambda T^*M)$ then $\alpha\in L^2$. And this implies
that an $H^1_0-$ harmonic form is in $L^2$.
We can generalize this fact with the Sobolev Orlicz inequality we have shown in \cite{carrcol}
\begin{theo}There is an constant $C(n)$ such that if $(M^n,g)$ is a complete Riemannian Manifold
whose curvature tensor $R$ and whose heat kernel $\left(P(t,x,y)\right)_{(t,x,y)\in \reel^*_+\times M\times M}$ satisfy
$$\int_M |R|^2(x)\left(\int^\infty_{C(n)/ |R|(x)}\sqrt{P(t,x,x)}dt\right)^2dx<\infty$$
then the exact sequence (\ref{eq:suitexa2}) holds.
\end{theo}
\proof First we remark that our assumption implies that we have
$$\int_M |R|(x)\left(\int^\infty_{C(n)/ |R|(x)}\sqrt{\frac{P(t,x,x)}{ t}}dt\right)^2dx<\infty$$
So that there is a choice of the constant $C(n)$ which implies that the Gauss-Bonnet operator is non-parabolic at infinity by the proposition
 (\ref{the:orlicznpy}). 

Now we have only to show that if $\alpha\in H^1_0$ is a harmonic form then $\alpha$ is in $L^2$,
 because we have in our case 
$H^1_0=W$. 
For this we used the Bochner-Weitzenb\"ock formula
$$\Delta\alpha=\bar \Delta\alpha+\cou \alpha=0.$$
The main point is the following :

{\it if $K$ is a compact of $M$, big enough and $\Omega=M\setminus K$, then we have
that
\begin{eqnarray*}&\|\bar\Delta_\Omega^{-1}\cou\|_{L^2(\Omega)\to L^2(\Omega)}<1\\
            &\|\bar\Delta_\Omega^{-1}\cou\|_{H^1_0(\Omega)\to H^1_0(\Omega) }<1.
            \end{eqnarray*}
Where $\bar \Delta_\Omega$ is the minimal self adjoint extension of
$\bar\Delta \,:\,C^\infty_0(\Lambda T^*\Omega)\la C^\infty_0(\Lambda T^*\Omega)$.}

For proving this we notice that we have the following 
$$|\cou(x)|\le c(n)|R(x)|$$
where $R$ is the Riemannian curvature tensor ; and by a consequence of Kato's lemma, we have
$$|\bar\Delta^{-1}_\Omega\alpha|(x)\le \Delta_\Omega^{-1}|\alpha|(x).$$
See for instance, the paper of \cite{HSU} or the appendix of \cite{Berard} written by G. Besson. So that we have only to show that
\begin{eqnarray*}\|\Delta_\Omega^{-1} |R|\|_{L^2(\Omega)\to L^2(\Omega)}<1/c(n)\\
            \|\Delta_\Omega^{-1}|R|\|_{H^1_0(\Omega)\to H^1_0(\Omega) }<1/c(n).
  \end{eqnarray*}          
Where $ \Delta_\Omega$ is the self adjoint extension of the essentially self adjoint operator
$\Delta \,:\,C^\infty_0(\Omega)\la C^\infty_0(\Omega)$. Let us show the first bound, we're going to prove that
$\||R|\Delta^{-1}_\Omega\|_{L^2\to L^2}<1/c(n)$ ; which is the same result because an operator and its adjoint have the same norm.
Recall that in \cite{carrcol}, we have shown the following Sobolev Orlicz inequality :
If $u\in C_0^\infty(M)$, we associate to it 
$$I(u)=\inf\{\int_M uv;\  v\in C^\infty_0(M)\  s.t.\  \int_M \phi(|v(x)|,x)dx\le 1\}$$
Where $\phi$ is the function on $\reel_+\times M$ defined by
$$\phi(\lambda,x)=\lambda\left (\int^\infty_{1/ \sqrt\lambda}\sqrt{P(t,x,x)}dt\right)^2.$$
Where we have assume that this integral is finite, otherwise all what we will say is empty. $I$ is a norm and
the completion of the space $C^\infty_0(M)$ with respect to $I$ is a Banach space (an Orlicz space), make of locally
integrable function. We have shown that for a universal constant $C$, which is taken to be the same
 as the one appearing in the Sobolev inequality
of proposition (\ref{the:orlicznpy}), we have the Sobolev inequality
$$\forall u\in C^\infty_0(M),\ \ I(u^2)\le C \|\Delta u\|.$$
By definition, we have the H\"older inequality:
$$\int_\Omega |R|^2u^2\le
 I(u^2)\inf\left\{\lambda, \int_\Omega\phi\left({|R|^2(x)/ \lambda},x\right)dx\le 1\right\}\ ;$$
so that if we choose $\Omega$ in order that
$$\int_\Omega\phi({4c(n)^2C|R|^2(x)},x)dx\le 1,$$
then we obtain that
$$\ \forall u\in C^\infty_0(\Omega),\ \ \||R|u\|_{L^2}\le \frac{1}{ 2c(n)}\|\Delta u\|_{L^2},$$
this achieve the proof of the first bound. For the second bound, we have only to show that
$$\|\Delta^{-1/2}_{\Omega}|R|\Delta^{-1/2}_{\Omega}\|_{L^2\to L^2}<{1/c(n)}$$
this because $\Delta^{-1/2}_{\Omega}$ is an isometry between $L^2(\Omega)$ and $H^1_0(\Omega)$.
For this we have only to show that $\||R|^{1/2}\Delta^{-1/2}_{\Omega}\|_{L^2\to L^2}<{1/ \sqrt{c(n)}}$.
Recall that we have the Sobolev inequality
 $$ \forall u\in C^\infty_0(M),\ \ N(u^2)\le C\|du\|_{L^2},$$
Now we used the H\"older inequality
$$\int_\Omega |R|u^2\le N(u^2)\inf\{\lambda, \int_\Omega\varphi({|R|(x)/ \lambda},x)dx\le 1\}\ ;$$
Now because $\varphi(\lambda,x)\le\phi(\lambda^2,x)$, under our assumption,  we have 
$$ \forall u\in C^\infty_0(\Omega),\ \ \int_\Omega |R|u^2\le \frac{1}{ 2Cc(n)} N(u^2).$$
This ends the proof of our two bounds.

We can now finish the proof. If $\rho$ is a smooth function which is
$1$ out of compact and with support in $\Omega$ then the form
$(\bar\Delta+\cou)\rho\alpha=\beta$ is a smooth form with compact support in $\Omega$. We remark that

{\it  $\bar \Delta^{-1}_\Omega\beta$ is in $H^1_0$ and in $L^2$. }

As the matter of fact, the linear form
$\phi\mapsto<\beta,\phi>_{L^2}$ is continuous on $H^1_0$,we have an $\gamma\in H^1_0$ so that
$$<\gamma,\bar\Delta \phi>_{L^2}=<\beta,\phi>_{L^2},\ \forall \phi\in C^\infty_0(\Lambda T^*\Omega),$$
so that $\gamma=\bar\Delta_\Omega^{-1} \beta\in H^1_0$. For proving that $\bar\Delta_\Omega^{-1} \beta$
is in $L^2$ we note that if $f$ is a smooth function with compact support which is $1$ on the support of $\beta$ then
$$\bar\Delta_\Omega^{-1} \beta=\bar\Delta_\Omega^{-1} f\beta.$$
But by what we have done before we have that $\phi\mapsto \bar\Delta_\Omega^{-1} f\phi$
 is a bounded operator of $L^2$, so
$\bar\Delta_\Omega^{-1} \beta$ is in $L^2$.

 Now we wrote that 
$$({\rm Id}_{H^1_0}+\bar\Delta_\Omega^{-1}\cou)\rho\alpha=\bar\Delta_\Omega^{-1} \beta.$$
But the operator $({\rm Id}+\bar\Delta_\Omega^{-1}\cou)$ is invertible on $H^1_0$ and on $L^2$ so that
$$\rho\alpha=({\rm Id}_{H^1_0}+\bar\Delta_\Omega^{-1}\cou)^{-1}\bar\Delta_\Omega^{-1} \beta$$
But we have
$$ ({\rm Id}_{H^1_0}+\bar\Delta_\Omega^{-1}\cou)^{-1}\bar\Delta_\Omega^{-1}\beta=
({\rm Id}_{L^2}+\bar\Delta_\Omega^{-1}\cou)^{-1}\bar\Delta_\Omega^{-1}\beta\,,$$
so that $\rho\alpha$ is squared integrable and $\alpha $ is in $L^2$.\hfill\endproof

 In \cite{carrcre}, we have given a cohomological interpretation of the dimension of
$$\frac
{\{\alpha\in W,\, d\alpha+\delta\alpha=0\}}
{\{\alpha\in L^2,\, d\alpha+\delta\alpha=0\}}.$$
 This dimension can be thought as the dimension of half bound states in Quantum Mechanics.
  A trivial cases when we know it's zero, is when 
zero is not in the essential spectrum of the Laplacian, or equivalently when the Gauss-Bonnet
 operator is Fredholm on its
domain. But in this case, the spaces of  $L^2$-cohomology and of reduced $L^2$-cohomology are the same, and
the long exact sequence (\ref{eq:suitexa2}) always hold for the $L^2$-cohomology. For instance :
\begin{theo}If $M$ is an even
dimensional locally symmetric space of finite volume and
negative curvature then
the exact sequence \ref{eq:suitexa2} holds.
\end{theo}

This is a consequence of the work of Borel and Casselman \cite{BC}: the Gauss-Bonnet operator
 is a Fredholm operator in this case.

With the calculus we have done in \cite{carrcre} we can now complement the proposition (\ref{the:wp}) 
\begin{theo}
If $(M,g)$ is a complete Riemannian Manifold of dimension $n$ and if there is a compact set $K$ of $M$ such that
$(M\setminus K,g)$ is isometric to the warped product $]0,\infty[\times\partial K,dr^2+f^2(r)g)$,
 then the Gauss-Bonnet operator is non-parabolic at infinity 
in the following two cases
$f(r)=ar$ for an $a>0$, or $\lim_{r\to \infty}f'(r)=0$.
\begin{enumerate}
\item Moreover let $k=[{n/2}]$, if $f(r)=ar$, 
and if the first eigenvalue $\lambda_0^k$ of the Laplace operator of acting on closed $k$
-differential forms of $\partial K$ 
satisfies $\lambda^k_0>\frac{7+(-1)^n}{ 8}a$ then the exact sequence (\ref{eq:suitexa2}) hold. Moreover when the Manifold is oriented and 
of even dimension then $$\chi_{L^2} (M)=\chi(M)-\sum_{j=0}^{k-2}(-1)^j{\rm b}_j(\partial K).$$

\item If $\lim_{r\to \infty}f'(r)=0$ and if we have $ \int^\infty f^t<\infty$ for all $t>1$ then the exact sequence
hold in even dimension in that case we have $\chi_{L^2} (M)=\chi(M)+\sum_{j=0}^{k-1}(-1)^j{\rm b}_j(\partial K)\ ;$ and 
in odd dimension $\dim M=2k+1$, it hold, provided that moreover the $k^{\rm th}$ Betti number of $\partial K$ vanish
i.e. ${\rm b}_k(\partial K)=0$.\end{enumerate}
\end{theo}
In fact, the second Gauss-Bonnet formula was already known by J. Br\"uning ([Br]). 

\subsection{Novikov-Shubin invariants and non-parabolicity at infinity.}
We give a new spectral condition which implies that the Gauss-Bonnet operator is non-parabolic at infinity.
\begin{pro}Assume that $(M,g)$ is a complete Riemannian Manifold such that 
there is an $\alpha>2$ and
a locally bounded function $C(x), \, x\in M$ with
$$\|e^{-t\Delta}(x,x)\|\le C(x)t^{-\alpha/2},\ \forall t\ge 1, x\in M$$
then the Gauss-Bonnet operator of  $(M,g)$ is non-parabolic at infinity and $h_\infty(M)=0$.
\end{pro}
\proof We claim that under these assumptions the integral 
$$G=\int_0^\infty (d+\delta)e^{-t\Delta} dt$$
defines a bounded operator from $L^2$ to $H^1_{loc}$ and moreover 
$$(d+\delta)G={\rm Id}_{L^2}, G(d+\delta)={\rm Id}_{C^\infty_0}.$$
The second equality implies that $(d+\delta)$ is non-parabolic at infinity and both that
$(d+\delta)$ has no $L^2$-kernel. Moreover if we define $W(E)$ as the completion
 of the space $C_0^\infty(E)$ with respect to the norm
$$\sigma\mapsto \|(d+\delta)\sigma\|_{L^2};$$
then the injection of $C_0^\infty(E)$ into $H^1_{loc}$ extends by continuity to an injection of $W(E)$ into $H^1_{loc}$ ;
 since
$(d+\delta)\,:\,W\la L^2$ is an isometry, so that $(d+\delta)$ has no $W$-kernel.

 We have to show the convergence: for this we cut the integral at $t=1$
$$G=\int_0^1 (d+\delta)e^{-t\Delta} dt+\int_1^\infty (d+\delta)e^{-t\Delta} dt.$$
The first integral converges because of the spectral theorem, the second is convergent
in the norm topology in the space of bounded operator from $L^2$ to $L^{\infty}_{loc}$.
As a matter of fact, the spectral theorem shows that
$$\|(d+\delta)e^{-\frac{t}{2}\Delta}\|_{L^2\to L^2}\le\sqrt{2/ t};$$
moreover if $f\in L^2(E)$, we have 
\begin{eqnarray*} \|(d+\delta)e^{-t\Delta}f(x)\| &=\|e^{-\frac{t}{2}\Delta}(d+\delta)e^{-\frac{t}{2}\Delta}f(x)\|\\
                               &\le \sqrt{\|e^{-t\Delta}(x,x)\|}\sqrt{\frac{2}{ t}}\|f\|\cr
                               &\le \sqrt{2C(x)}t^{-\frac{1}{ 2}-\frac{\alpha}{ 4} }\|f\|.
  \end{eqnarray*}                             
And this proves the convergence of the second integral. The two equalities
 are consequences of our hypothesis and of the formula
$$(d+\delta)\int_0^T (d+\delta)e^{-t\Delta} dt={\rm Id}-e^{-T\Delta}.$$
 \hfill\endproof
 
We make some remarks :
\begin{itemize}
\item The  non-parabolicity at infinity property and the dimension $h_\infty(M)$
 depends only on the geometry at infinity, so
 the conclusion of the theorem remains valid for any complete Riemannian Manifold which is isometric 
 outside some compact set to one satisfying the hypothesis of the theorem.
\item By the Karamata theorem, our assumption is equivalent to the
 following on $E$ the spectral resolution of $\Delta$
$$\|E([0,\lambda],x,x)\|\le \tilde C(x) \lambda^{\alpha/ 2}, \forall \lambda \in [0,1], \forall x\in M.$$
\end{itemize}
In fact the best possible exponent $\alpha$ is linked with the Novikov-Shubin invariants :
if $(M,g)$ is the universal covering of a compact Manifold $\overline M$ and if $D=d+\delta$ is the Gauss-Bonnet
operator acting on differential forms on $M$, and if $F\subset M$
 is a fundamental domain of the covering $M\to {\overline{M}}$, then
$$\beta=\inf\{\alpha\,|\,\int_F{\rm trace}_{\Lambda T^*_x M}e^{-t\Delta}(x,x)dx=O(t^{-\alpha/2})\}$$
is a Novikov-Shubin invariant of $\overline M$,
 it doesn't depend on the metric (\cite{ns}), nor on the differential structure
(\cite{Lottns}) and it is a homotopy invariant of $\overline M$ (\cite{GS}).

According to the calculus done by M. Rumin (\cite{rumin}) and by L. Schubert (\cite{schu}) we have :
\begin{cor}
If $n>1$ then the Gauss-Bonnet operator of the Heisenberg group 
$$H_{2n+1}=\left\{
 \left(  \begin{array}{cccccc}
                    1&x_1&x_2&...&x_n&z\\
                    0&1  &0&...&0&y_1\\
                    .&.. &.&...&.&.\\
                    .&.. &...&.&.&.\\
                    .&.. &...&0&1&y_n\\
                    0&.. &.&...&.&1\\
                    \end{array}\right) ;\ \  x_1,..,x_n,y_1,..,y_n,z\in \reel\right \}$$ with a left invariant metric is non-parabolic at infinity .
                    Furthermore, for any complete Riemannian Manifold which is isometric outside 
                    some compact to such $H_{2n+1}$, the exact sequence (\ref{eq:suitexa2}) holds.
 \end{cor}                   
\proof It is shown in \cite{rumin} and \cite{schu} that
the hypothesis of our proposition (3.10) holds with $\alpha=n+1$ for the Gauss-Bonnet operator of 
the Heisenberg group $H_{2n+1}$ with a left invariant metric.\hfill\endproof

\section{ Manifolds with flat ends}

In this part, we apply our results to complete Riemannian Manifold with flat ends.
 Let $(M^n,g)$ be a complete Riemannian Manifold
whose curvatures vanish outside
 some compact set.
 \subsection{Geometry of flat ends.}
 Then such a Manifold has a finite number of ends $E_1,\, E_2,..., E_b$. And according to the works
 of Eschenburg-Schroeder such these flat ends are classified in three families \cite{ES}:
 \begin{theo}
 \begin{enumerate}
 \item $E$ has a finite cover isometric to $\left(\reel^{\nu}-B(R)\right)\times \T^{n-\nu}$, where $\nu>2$ and
 $\T^{n-\nu}$ is a flat $(n-\nu)$-torus.
 \item $E$ is isometric to the Riemannian product $]0,\infty[\times Y$ where $Y$ is a compact flat Manifold.
 \item The universal cover of $E$ is isometric to $Y_{\beta,R}\times \reel^{n-2}$, for $R>\beta$.
  And the $\pi_1(E)$ respects this decomposition.
 Where $Y_{\beta,R}$ is defined as follows:
 
  $c_{\beta,R}(s)=(\beta s-R \sin s, R\cos s)$ is the cycloid and $n_{\beta,R}$ its
 unit normal vector. Then $Y_{\beta,R}$ is $]0,\infty[\times\reel$ equipped with the metric $\Phi_{\beta,R}^*(dx^2+dy^2)$
 where $\Phi_{\beta,R}\, :\, ]0,\infty[\times\reel\longrightarrow \reel^2$
 is the immersion $\Phi_{\beta,R}=c_{\beta,R}(s)+t n_{\beta,R}(s)$.
 \end{enumerate}
 \label{classi} \end{theo}
 In the third case, a finite cover of the end is isometric to
   $\left(Y_{\beta,R}\times \reel^{n-2} \right)/ \Gamma$ a flat bundle of $(n-2)$ flat torus on
 $\bar Y_{\beta,R}=Y_{\beta,R}/\{(s,t)\sim (s+k \tau,t), k\in \ent\}$ with $\tau=2\pi$ if $\beta>0$ and 
 $\tau$ is any positive real number if $\beta=0$.
   
\subsection{The $L^2$-cohomology of flat ends.}
We want to apply our long exact sequence to describe the $L^2$-cohomology of Manifolds
with flat ends, so we begin to compute the $L^2$-cohomology of flat ends.

\begin{lem}
In case (2) of the classification (\ref{classi}), 
if $\#$ is the relative or absolute boundary condition, then $H^k_{2,\#}(E)=\{0\}$.
\label{the:nulcyl}
\end{lem} 
\proof
According to the proposition (\ref{the:double}), it is enough to show the $L^2$ cohomology of the double Manifold
$E\#E$ has a trivial $L^2$-cohomology. But $E\#E$ is the Riemannian product $\reel\times Y$. And it is well known that 
$H^k_{2}(\reel\times Y)=\{0\}$ (this is for example a consequence of the Kunneth formula).
\endproof

This lemma can also be recover from the result of Atiyah-Patodi-Singer (\cite{APS})
which that on a Manifold with cylindrical end, the $L^2$-cohomology is the image of
the relative cohomology in the absolute cohomology.
\begin{lem}
In case (3) of the classification (\ref{classi}), 
if $\#$ is the relative or absolute boundary condition, then $H^k_{2,\#}(E)=\{0\}$.
\label{the:nulqua}
\end{lem} 
\proof
Let $\pi\,:\,\widehat{ E}\longrightarrow E$ the finite cover of $E$ which is a flat bundle of $(n-2)$ flat torus on
 $\bar Y_{\beta,R}$. We have 
 $$
\alpha\in L^2(\Lambda T^*E) {\rm\ if\ and\ only\ if\ } \pi^*\alpha \in L^2(\Lambda T^*\widehat{E})$$
So it is enough to show the result for $\widehat{ E}$.
Now, the $(n-2)$ flat torus acts on $\widehat{ E}$ by isometry and this action generate Killing fields of bounded length, 
so by the argument of
N. Hitchin, a $\#$-$L^2$ harmonic forms must be invariant by the torus action (theorem 3 of \cite{hitchin}).
So that there is an identification between $H^k_{2,\#}(\widehat{ E})$ and the space of $L^2$ harmonic forms on 
$\bar Y_{\beta,R}$ with value in the flat unitary bundle
 $\left(Y_{\beta,R}\times \Lambda (\reel^{n-2})^* \right)/ \Gamma$ satisfying the relative or absolute boundary condition.
 Now there is no non trivial such $0$ and $2$ forms as they must be parallel. And for $1$-forms being $L^2$ harmonic is a conformal
 invariant property. As $\bar Y_{\beta,R}$ is conformally equivalent to the cylinder, the lemma (\ref{the:nulcyl})
 implies that there is no $L^2$ harmonic $1$- forms with value in this flat unitary bundle.
\endproof
\begin{lem}
In case (1) of the classification (\ref{classi}), if $E$ is isometric to the product $\left(\reel^{\nu}-B(R)\right)\times \T^{n-\nu}$
 then $$H^k_{2,rel}(E)=\reel \frac{dr}{r^{\nu-1}} \wedge H^{k-1}(\T^{n-\nu}) .$$
 \label{the:nphodge}
\end{lem} 
\proof
The proposition 4.3 of \cite{carrma} implies that 
$H^k_{2,rel}(\reel^{\nu}-B(R))=\reel \frac{dr}{r^{\nu-1}}$;
hence the lemma by the Kunneth formula.
\endproof
\subsection{The parabolic ends} We begin by
eliminated the parabolic ends. Let 
$E_P$ the union of parabolic ends of $M$, i.e. the union of the ends of
type (2) and (3) in the classification (\ref{classi}).

\begin{pro}
$$H^k_{2}(M)\simeq \ima\left(  H^k_{2,rel}(M\setminus E_P) \longrightarrow H^k_{2,abs}(M\setminus E_P)\right).$$
\end{pro}
\proof
First if we look at the short sequence 
$$ H^k_{2,rel}(M\setminus E_P) \stackrel{\mbox{i}} {\longrightarrow}
 H^k_{2}(M)\stackrel{\mbox{j}^*}{\longrightarrow} H^k_{2,abs}(E_P).$$
From the theorem (\ref{the:seb}), we know that $\ker \mbox{j}^*=\ima\mbox{i}$ and from
the lemmas (\ref{the:nulcyl},\ref{the:nulqua}), we know that $H^k_{2,abs}(E_P)=\{0\}$. 
So in any $L^2$ cohomology class on $M$ there is a closed $L^2$ smooth form with support in $M\setminus E_P$.
Second, we look at the short sequence
$$ H^k_{2,rel}(E_P) \stackrel{\mbox{i}} {\longrightarrow}
 H^k_{2}(M)\stackrel{\mbox{j}^*}{\longrightarrow} H^k_{2,abs}(M\setminus E_P).$$
Again, this sequence is exact and $H^k_{2,rel}(E_P)=\{0\}$.

So the injective map between $\ima\left( H^k_{2,rel}(M\setminus E_P) \rightarrow H^k_{2,abs}(M\setminus E_P)\right)$
and $H^k_{2}(M)$ is surjective hence it is an isomorphism. \endproof

\subsection{The non-parabolic ends} We have now to compute the reduced 
$L^2$ cohomology of $M\setminus E_P$ with either relative or absolute boundary condition on
$\partial E_P$. We'll treat only the case of absolute boundary conditions ; the case of
relative boundary conditions is the same.

{\it  Now on $(M\setminus E_P,g)$ any $W$-harmonic form must be $L^2$:}

As a matter of fact, by (\ref{the:eind}) this property only depends on the geometry of ends: if
$E$ is an end of $M\setminus E_P$, then $E$ has a finite cover isometric to
$\widehat{E}=\left(\reel^{\nu}-B(R)\right)\times \T_E$, where $\nu>2$ and
 $\T_E$ is a flat $(n-\nu)$-torus. Now a harmonic form on $E$ is in $W$ (resp. $L^2$) if and only if
 it is pulled back to a
 $W$ (resp. $L^2$) form on $\widehat{E}$ . And $\widehat{E}$ is the end of $\reel^{\nu}\times \T_E$ so by the theorem
 (\ref{the:eind}), it is enough to show the result for $\reel^{\nu}\times \T_E$. But
 in (\cite{carrpac}) it is shown that
 $\ker_W (d+\delta)_{\reel^{\nu}\times \T_E}=\ker_W (d+\delta)_{\reel^{\nu}}
 \otimes \ker (d+\delta)_{ \T_E}$.
 
Now, for $\nu>2$, $\reel^{\nu}$ is non-parabolic so the topology of $W$ is given by the quadratic forms
$\alpha\mapsto \int_{\reel^{\nu}}|\nabla\alpha|^2=\int_{\reel^{\nu}}|(d+\delta)\alpha|^2$; so a $W$-harmonic form is
parallel. But by the Sobolev inequality, an element of $W$ has to be in $L^{2\nu/(\nu-2)}$. Hence any $W$-harmonic form
 on $\reel^{\nu}$ is zero.
 
So on $M\setminus E_P$, we have $h_\infty(M\setminus E_P)=0$ and let $E_{NP}=\cup_{i=1}^b E_i$ be the union of the ends of $M\setminus E_P$ and
$U=E_P\cup E_{NP}$ be the union of the ends of $M$.
We have the long exact sequence :

$$..\la H^k_{2,rel}(E_{NP}) \stackrel{\mbox{i}}{\longrightarrow} H^k_{2,abs}(M\setminus E_P)
 \stackrel{\mbox{j}^*}{\longrightarrow} H^k_{abs}(M\setminus U)
\stackrel{\mbox{b}}{\longrightarrow} H^{k+1}_{2,rel}(E_{NP})\la..$$

Each $E_i$ has a finite cover isometric to
$\widehat{E_i}=\left(\reel^{{\nu_i}}-B(R)\right)\times \T_i$ where $\nu_i>2$ and $\T_i$ is a flat $(n-{\nu_i})$-torus.
So for a finite subgroup $\Gamma_i$ of $O({\nu_i})\times {\rm Isom} (\T_i)$, 
$E_i=\widehat{E_i}/\Gamma_i $. Let $p_i\,:\, \widehat{E_i}\longrightarrow E_i$
be the covering map. $p_i^*$ induced an isomorphism between  $H^k_{2,rel}( E_i)$
 and the space of $\Gamma_i$-invariant element in $H^k_{2,rel}(\widehat{E_i})$:
$$H^k_{2,rel}( E_i)=H^k_{2,rel}(\widehat{E_i})^{\Gamma_i}.$$  We have shown that
$$H^k_{2,rel}(\widehat{E_i})\simeq \frac{dr}{r^{{\nu_i}-1}}\wedge H^{k-1}(\T_i).$$
And the isomorphism between $H^k_{2,rel}(\widehat{E_i})$ and $H^{k-1}(\T_i)$ is induced by the map
$$ H^{k-1}(\T_i)\rightarrow H^{k-1}(\partial \widehat{E_i})\stackrel{\mbox{b}}{\longrightarrow} H^k_{2,rel}(\widehat{E_i}).$$
$\Gamma_i$ acts on $\T_i$ so $$H^k_{2,rel}( E_i)=H^{k-1}(\T_i)^{\Gamma_i}.$$
Let $\pi\,:\, \cup_i\T_i\longrightarrow \partial E_{NP}$ be the immersion induced by
the $p_i$'s. And define
$\Omega\left( M\setminus U, \ker \pi^*\right)$ be the space of smooth differential form on
$ M\setminus U$ which are zero when pull back by $\pi$ to each $\T_i$. Then
$\Omega\left(M\setminus U, \ker \pi^*\right)$ is a subcomplex of the complex of differential forms
on $ M\setminus U$. Let $H^k\left(M\setminus U,\ker \pi^*\right)$
the associated cohomology spaces. Now we have also the exact sequence

$$..\la \bigoplus_i  H^{k}(\T_i)^{\Gamma_i}\stackrel{\mbox{i}}{\longrightarrow} H^k\left(M\setminus U, \ker \pi^*\right)
 \stackrel{\mbox{j}^*}{\longrightarrow} H^k_{abs}(M\setminus U)
\stackrel{\mbox{b}}{\longrightarrow} \bigoplus_i  H^{k+1}(\T_i)^{\Gamma_i}\la..$$

We can build a map $H^k_{2,abs}(M\setminus E_P)\la H^k\left(M\setminus U, \ker \pi^*\right)$:
if $\alpha$ is a smooth closed $L^2$ form on $M\setminus E_P$, then from (\ref{pro:greenb}) and (\ref{the:nphodge})
then on the end $E_i\subset E_{NP}$, we have
 $$\alpha=\frac{dr}{r^{{\nu_i}-1}}\wedge \beta_i+d\gamma_i$$
 where $\gamma\in W(\Lambda^{k-1}T^*E_{i})$ and $r$ is the radial coordinate on $E_{i}$, on each end it is 
 the function distance to $\partial E_{i}$ plus a constant.
 Now let $\overline{\gamma}\in W(\Lambda^{k-1}T^*(M\setminus E_P))$ be an extension of the $\gamma_i$'s, then
 $\alpha-d\overline{\gamma}$ is a closed element of $\Omega\left( M\setminus U, \ker \pi^*\right)$ and its
 cohomology class in $H^k\left(M\setminus U,\ker \pi^*\right)$ depends only on the $L^2$
 cohomology class of $\alpha$. So we get the desired map.
 Now, we have a commutative diagram
$$
\xymatrix{
  {H^{k-1}_{abs}(M\setminus U)} \ar[d] \ar[r]^{b} & {H^k_{2,rel}(E_{NP}) }\ar[d]\ar[r]^{i}
  & {H^k_{2,abs}(M\setminus E_P)}\ar[d] \ar[r]^{j^*} & {H^k_{abs}(M\setminus U )} \ar[d]\ar[r]^{b}     & {H^{k+1}_{2,rel}(E_{NP})} \ar[d] \\
  {H^{k-1}_{abs}(M\setminus U)} \ar[r]        & {\bigoplus_i  H^{k-1}(\T_i)^{\Gamma_i}} \ar[r]^{b}     &
  {H^k\left(M\setminus U, \ker \pi^*\right)}\ar[r]^{j^*}      & {H^k_{abs}(M\setminus U )}\ar[r] &
  {\bigoplus_i  H^{k+1}(\T_i)^{\Gamma_i}}
  }
$$

 The first two and last two vertical arrows are isomorphisms, hence by the fifth arrow lemma, we have an isomorphism between
 $H^k_{2,abs}(M\setminus E_P)$ and  $H^k\left(M\setminus U, \ker \pi^*\right)$. We have prove :
 
\begin{pro}
$$H^k_{2,abs}(M\setminus E_P)\simeq H^k\left(M\setminus U, \ker \pi^*\right).$$
\end{pro}

And we arrive to our main result:

\begin{theo} Let $(M^n,g)$ be a complete Riemannian Manifold with flat ends. 
$M$ has a finite number of ends, let $E_P$ the union of the ends with sub-quadratic
 growth and $E_{NP}$ be the union of the remaining ends.
Then each connected component of $\partial E_{NP}$ has a finite cover isometric to
$\Sphere^{\nu-1}(R)\times \T^{n-\nu}$ where $\nu>2$ and $\T^{n-\nu}$ is a flat torus.
Let $\Omega(M\setminus (E_P\cup E_{NP}), A)$ the sub-complex of the complex of differential forms on
$M\setminus (E_P\cup E_{NP})$, consisting of forms which are zero when pulled back to these torii and
 let $\Omega(M\setminus (E_P\cup E_{NP}), A_c)$ be the subcomplex of this complex consisting of forms
which are zero when pulled back to $\partial E_P$. We denote by $H^k((M\setminus (E_P\cup E_{NP}), A)$ and
$H^k((M\setminus (E_P\cup E_{NP}), A_c)$ the associated cohomology spaces. Then we have the isomorphism :
$$H^k_{2}(M)\simeq \ima \left( H^k((M\setminus (E_P\cup E_{NP}), A_c)\rightarrow H^k((M\setminus (E_P\cup E_{NP}), A)\right)$$

\end{theo}

\subsection{A $L^2$ Chern-Gauss-Bonnet formula}
 
 We investigate now a Chern-Gauss-Bonnet formula for the $L^2$ Euler characteristic.
 We consider Manifold with only non-parabolic flats ends. In this case, the exact sequence (\ref{eq:suitexa2}) hold, hence
 we have
  $$\chi_{L^2}(M)=\chi(M)+\chi_{L^2,rel}(E_{NP}).$$
$E_{NP}$ is the union of the ends of $M$, each of these ends $E$ has a finite cover isometric to
$\widehat{E}=\left(\reel^{\nu}-B(R)\right)\times \T_E$ where $\T_E$ is a flat $(n-\nu)$ torus.
And for a finite subgroup $\Gamma_E$ of $O(\nu)\times {\rm Isom} (\T_E)$, 
$E=\widehat{E}/\Gamma_E $.
Assume that $M$ is oriented and  even dimensional then we can compute $\chi (M)$ with the Chern formula :
$$\chi(M)=\int_M\Omega^g+\sum_E\lim_{ \rho\to \infty} \int_{\Sigma_E(\rho)} P(II);$$
where $\Omega$ is the Euler form of $(M,g)$ and the summation is with respect to all
ends  $E$ of $M$,  $\Sigma_E(\rho)$ is the hypersurface $ (\rho \Sphere^{\nu-1}\times \T_E)/\Gamma_E$ and $P(II)$ is 
a polynomial expression of the curvature and of the second fundamental form of 
$\Sigma_E(\rho)$. We can compute $\int_{\Sigma_E(\rho)} P(II)$. This integral is multiplicative
 with respect to finite cover so that we have :
$$\int_{\Sigma_E(\rho)} P(II)=\frac{1}{\#\Gamma_E}\int_{\rho \Sphere^{\nu-1}\times \T_E} P(II).$$
But again by Chern formula, we have
$$\int_{\rho \Sphere^{\nu-1}\times \T_E} P(II)=\chi(\rho \Disk^{\nu-1}\times \T_E)=\chi( \Disk^{\nu})\, \chi(\T_E) ,$$
 so that we conclude
$$\int_{\Sigma_(E)(\rho)} P(II)=
\left\{
 \begin{array}{lcl}
0 & {\rm if } &\nu< \dim M \\
\frac{1}{\#\Gamma}&{\rm if } &\nu=\dim M.\\

  \end{array}    
  \right.$$

We compute now the $L^2$ Euler characteristic of an end $E$ of $M$:
we know that $$H^k_{2,rel}(\widehat{E})=\frac{dr}{r^{\nu-1}}\wedge H^{k-1}(\T_E).$$
Let $G_E$ be the image of $\Gamma_E$ in ${\rm Isom} (\T_E)$ then we have also that
$$H^k_{2,rel}(\widehat{E})\simeq H^{k-1}(\T_E)^{G_E}.$$
And we get $\chi_{L^2}(E)=-\chi(\T_E,G_E)$, where $\chi(\T_E,G_E)$ is the $G_E$ equivariant Euler characteristic of $\T_E$.
We can give a more precise formula. We know that the cohomology of the torus $\T_E$ is the
exterior algebra of $H^1(\T_E)$ or the exterior algebra of left invariant differential forms on $\T_E$.
But we have $$\dim H^{k}(\T_E)^{G_E}=\frac{1}{\# G_E}\sum_{\gamma\in G_E} {\rm Tr}_{\H^k(\T_E)}\gamma^*$$
so $$\chi(\T_E,G_E)=\sum_k (-1)^k \dim H^{k}(\T_E)^{G_E}=
\frac{1}{\# G_E}\sum_{\gamma\in G_E}\sum_k (-1)^k{\rm Tr}_{\H^k(\T_E)}\gamma^*.$$
We have the formula :
$$\sum_k (-1)^k{\rm Tr}_{\H^k(\T_E)}\gamma^*=\det (Id_{\H^1(\T_E)}-\gamma^*).$$
This formula can also be recover from the Lepschetz fixed point formula.  
And we have prove the following theorem
\begin{theo} If $(M^n,g)$ is a complete oriented Riemannian Manifold of even dimension whose
curvatures vanish outside some compact and if for each end $E$ of $M$ we have
$$\lim_{r\to\infty}\frac{\vol E\cap B_x(r) }{r^2}=\infty,$$ then
$$\chi_{L^2}(M)=\int_M \Omega^g+
\sum_{E{\rm\ end\ of\ } M}q(E),$$
where $q(E)$ is defined as follow :
\begin{itemize}
\item When $\pi_1(E)$ has no torsion then $q(E)=0$.
\item  When ${\rm rank}\, \pi_1(E)=0$ we have $q(E)=\frac{1}{ |\pi_1(E)|}-1$
\item  When ${\rm rank}\,\pi_1(E)>0$ then $\pi_1(E)$ acts isometrically on $\Sphere^{\nu-1}\times \reel^{n-\nu}$, 
$n-\nu={\rm rank}\,\pi_1(E)<n-1$.  We can consider that $\pi_1(E)$ is a subgroup of
 $\ \Isom(\Sphere^{\nu-1}\times \reel^{n-\nu})=O(\nu)\times \left [\reel^{n-\nu}\rtimes O(n-\nu)\right]$ and let ${G_E}$ 
be the image of $\pi_1(E)$ in $O(n-\nu)$, then we have
$$q(E)=-\frac{1}{ |{G_E}|}\sum_{\gamma\in {G_E}}\det({\rm Id}-\gamma)$$ 
\end{itemize}
\end{theo}
When the Manifold has a parabolic end, the exact sequence doesn't hold and we cannot give an explicit formula
for the $L^2$ Euler characteristic. As a matter of fact, we cannot expect that the 
 $L^2$ Euler characteristic is the sum of the integral of the Euler form and of a contribution of ends.
A counterexample is given in (\cite{carrpac}): If $M=\reel^2\#\reel^2$ is the surface obtained by gluing two copies
of the Euclidean plane $\reel^2$; this surface has two planar ends and no non trivial  $L^2$-harmonic forms.
On $\reel^2$, we have $\chi_{L^2}(\reel^2)=0=\int_{\reel^2} \frac{KdA}{2\pi}$.
But on $M$, we have $\chi_{L^2}(M)=0\not = -2=\int_{M} \frac{KdA}{2\pi}$. The surfaces $M$ and $\reel^2\cup \reel^2$
 have the same ends, but the difference between the $L^2$ Euler characteristic and 
of the integral of the Euler form is not the same.

\end{document}